\documentclass[12pt]{amsart}
\usepackage{psfrag,epsfig,amssymb,url}
\numberwithin{figure}{section}
\numberwithin{table}{section}

\newtheorem{theorem}{Theorem}[section]

\newtheorem{lemma}[theorem]{Lemma}

\newtheorem{cor}[theorem]{Corollary}
\theoremstyle{definition}
\newtheorem{definition}[theorem]{Definition}

\theoremstyle{remark}
\newtheorem{remark}[theorem]{Remark}
\numberwithin{equation}{section}

\def \h{{\mathfrak h}}
\def \Z{{\mathbb Z}}

\def \[{[ }
\def \]{] }

\def \L{{\mathcal L}}

\begin{document}

\author{Anna Felikson}
\address{Independent University of Moscow}
\curraddr{School of Engineering and Science, Jacobs University Bremen}
\email{felikson@mccme.ru}
\thanks{Research of the second author is supported by grants RFBR 10-01-00678 and NSh 709.2008.1}

\author{Sergey Natanzon}
\address{Independent University of Moscow,  \phantom{wwwwwwwwwwwwwwwwwwwwwwwwwwww}
Belozersky Institute of Physico-Chemical Biology, Moscow State University,\phantom{wwwww}
Institute Theoretical and Experimental Physics}
\email{natanzons@mail.ru}

\title{Labeled double pants decompositions}

\begin{abstract} 
A double pants decomposition of a 2-dimensional surface 
is a collection of two pants decomposition of this surface introduced in~\cite{FN}.  
There are two natural operations acting on double pants decompositions: flips and handle-twists. It is shown in~\cite{FN} 
that the groupoid generated by flips and handle-twists acts transitively on admissible double pants  decompositions
where the class of admissible decompositions has a natural topological and combinatorial description.
In this paper, we label the curves of double pants decompositions and show that for all but one surfaces
the same groupoid acts transitively on all labeled admissible double pants decompositions. 
The only exclusion is a sphere with two handles, where the groupoid has 15 orbits.

\end{abstract}

\maketitle


\section*{Introduction}
Consider a 2-dimensional orientable surface $S$ of genus $g$ with $n$ holes. 
A pants decomposition of $S$ is a decomposition into 3-holed spheres  
(called ``pairs of pants''). 
In~\cite{FN} we considered double pants decompositions of surfaces as a union of two pants decompositions (with an additional assumption
that the homology classes of the curves contained in the double pants decomposition generate the whole homology lattice $H_1(S,\Z)$).
We introduced a simple groupoid acting on double pants decompositions (the groupoid is generated by transformations of two
types called  flips and handle-twists, each flip or handle-twist affecting only one curve of double pants decomposition)  
and proved that this groupoid acts transitively on all 
admissible double pants decompositions. The class of admissible double pants decompositions has a simple combinatorial definition
(see  Definition~\ref{admissible} below) as well as a nice description in terms of Heegaard splittings of 3-manifolds.   

More precisely, for each pants decomposition $P$ of $S$ one may construct a handlebody $S_+$ such that $S$ is 
the boundary of $S_+$  and all curves of $P$ are contractible inside $S_+$. A union of two pants 
decompositions of the same surface define two different handlebodies bounded by $S$, attaching this
handlebodies along $S$  one obtains a Heegaard splitting of some 3-manifold.
This connection of two pants decompositions to a Heegaard splitting was investigated in a row of papers (\cite{CG},~\cite{He},~\cite{L})
and many others, see~\cite{C} for further references). 
The two pants decompositions are considered usually as two vertices in a pants complex, using as the main tool the Hempel distance.

The admissible double pants decomposition defined in~\cite{FN} are exactly ones resulting in Heegaard splittings of
a 3-sphere. So, the transitive action of flip and handle-twists groupoid on admissible double pants decompositions may be interpreted as 
an action on Heegaard splitting of 3-sphere.

\medskip

In this paper we consider double pants decompositions with curves labeled by distinct  integer numbers. We define a trivial action of 
flips and handle-twists on the labels: all labels are preserved by these transformations, in particular, the label of the flipped or
twisted curve coincides with its initial label. We consider the action of the groupoid generated by flips and handle-twist on
labeled admissible double pants decompositions and obtain the following theorem:

\medskip
\noindent
{\bf Theorem A. (Main Theorem)}
{\it
The flip-twist groupoid acts transitively on labeled admissible double pants decompositions of $S_{g,n}$, $2g+n\ge 2$, unless $(g,n)=(2,0)$.
The action of flip and twist groupoid on labeled admissible double pants decompositions of $S_{2,0}$ has  15 orbits.
}
\medskip

Furthermore, we also may restrict ourselves to the case of one labeled pants decomposition. 
It is shown by Hatcher and Thurston~\cite{HT},~\cite{H1} that there are two types of transformations called 
flips and $S$-moves which are sufficient to connect all pants decompositions in the unlabeled case.  
We extend the statement to the labeled case:

\medskip
\noindent
{\bf Theorem B.}
{\it
The groupoid generated by flips and $S$-moves 
acts transitively on labeled pants decompositions of $S_{g,n}$, $2g+n\ge 2$.
}
\medskip

\medskip

The paper is organized as follows. In Section~\ref{dp} we recall from~\cite{FN} the definitions concerning double pants decompositions
and flip-twist groupoid. We also introduce the notion of labeled double pants decomposition and define the action of flip-twist groupoid 
on the labels. For the aims of proofs  we consider also a notion of strict-labeled double pants decompositions for which the action of 
flip-twist groupoid is unable to intermix the labels of one pants decomposition with the labels of another.
In Section~\ref{sec ordinary} we consider labeled pants decomposition and prove Theorem~B.
In Section~\ref{sec strict} we prove transitivity of flip-twist groupoid on strict-labeled decompositions
(to be exact, the groupoid acts transitively unless $(g,n)=(2,0)$ and has 6 orbits otherwise).
Finally, in Section~\ref{sec weak} we use the result of  Section~\ref{sec strict} to prove the Main Theorem.

\medskip
\noindent
{\bf Acknowledgments.} We are grateful to Robert Penner for suggestion to consider the orbits of labeled double pants decompositions.

\section{Double pants decompositions}
\label{dp}
In this section we introduce double pants decompositions and their transformations.

\subsection{Pants decompositions}

Let $S=S_{g,n}$ be an oriented closed surface of genus $g\ge 0$ with $n$ holes.
A {\it curve} $c$ on $S$ is an embedded closed non-contractible curve considered up to a homotopy of $S$. 

Given a set of curves we always assume that there are no ``unnecessary intersections'', so that if two curves of this set
intersect each other in $k$ points then there are no homotopy equivalent pair of curves intersecting in less than $k$ points.

For a pair of curves $c_1$ and $c_2$ we denote by $|c_1\cap c_2|$ the number of (geometric) intersections of $c_1$ with $c_2$.

\begin{definition}[{{\it Pants decomposition}}]
A {\it pants decomposition} of $S$ is a system of (non-oriented) mutually disjoined curves 
$P=\langle c_1,\dots,c_m \rangle$  decomposing  $S$ into pairs of pants 
(i.e. into spheres with 3 holes).

\end{definition}

It is easy to see that any pants decomposition of a surface $S_{g,n}$ consists of $m=3g-3+n$ curves.
For abuse of notation we will always write $m$ instead of $3g-3+n$.

Note, that we do allow self-folded pants, two of whose boundary components are identified in $S$.  
A surface which consists of one self-folded pair of pants will be called a {\it handle}.

\begin{definition}[{{\it Lagrangian plane of pants decomposition}}]
Let $P=\langle c_1,\dots,c_{m}\rangle$ be a pants decomposition. A Lagrangian plane  $\L(P)\subset H_1(S,\Z)$ is a subspace spanned 
by the homology classes $h(c_i)$, $i=1,\dots,m$ (here $c_i$ is taken with any orientation).

\end{definition}


\begin{definition}[{{\it Flip}}]
Let $P=\langle c_1,\dots,c_m  \rangle$ be a pants decomposition.
Define a {\it flip  of $P$ in the curve} $c_i$  as  
a replacing of a regular curve $c_i\subset P$ by any curve $c_i'$ satisfying the following properties:
\begin{itemize}
\item $c_i'$ does not coincide with any of $c_1,\dots,c_m$;
\item $|c_i'\cap c_i|=2$;
\item $c_i'\cap c_j=\emptyset$ for all $j\ne i$.
\end{itemize}

\end{definition}

\begin{figure}[!h]
\begin{center}
\psfrag{c1}{\scriptsize $c'$}
\psfrag{c}{\scriptsize $c$}
\epsfig{file=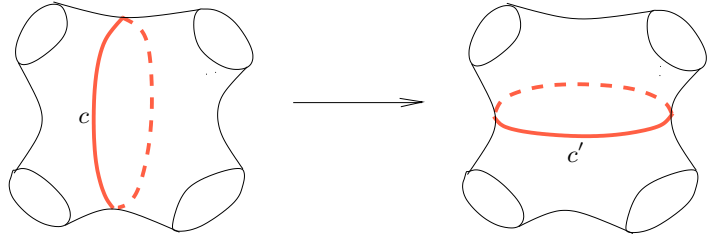,width=0.6\linewidth}
\caption{Flips of pants decomposition.} 
\label{fig-flip}
\end{center}
\end{figure}

See Fig.~\ref{fig-flip} for an example of a flip. Clearly, an inverse operation to a flip is also a flip (so that the set of flips
compose a groupoid acting on pants decompositions).

\subsection{Double pants decompositions }

\begin{definition}[{{\it Lagrangian planes in general position}}]
Two {\it  Lagrangian planes $\L_1$ and $\L_2$  are in general position} if 
 $\L_1\cap \L_2=0$ and $H_1(S,\Z)=\langle \L_1,\L_2\rangle$).

\end{definition}

\begin{figure}[!h]
\begin{center}
\psfrag{P}{$P_a$}
\psfrag{Z}{$P_b$}
\epsfig{file=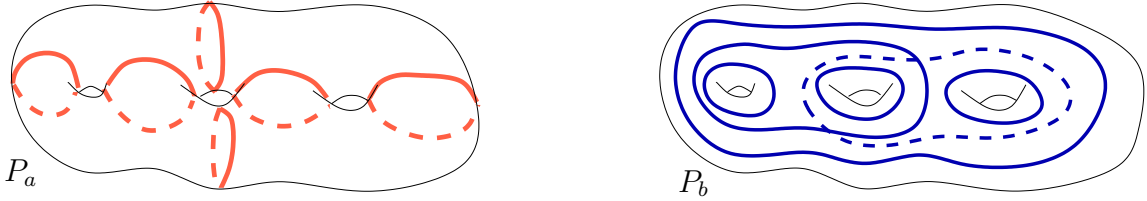,width=0.98\linewidth}
\caption{Pair of pants decompositions $(P_a,P_b)$. } 
\label{double_pant}
\end{center}
\end{figure}

See Fig.~\ref{double_pant} for an example of two pants decompositions spanning a pair of Lagrangian planes in general position.

\begin{definition}[{{\it Double pants decomposition}}]
A {\it double pants decomposition} $DP=(P_a,P_b)$ is a pair of pants decompositions $P_a$ and $P_b$ of the same surface
such that the Lagrangian planes $\L_a=\L(P_a)$ and  $\L_b=\L(P_b)$ spanned by these pants decompositions
are in general position.

\end{definition}

There are several natural transformations on the set of double pants decompositions:
\begin{itemize}
\item flips of $P_a$;
\item flips of $P_b$;
\item handle-twists (see Definition~\ref{Dehn} below). 

\end{itemize}

\begin{figure}[!h]
\begin{center}
\psfrag{a}{\scriptsize $a$}
\psfrag{b}{\scriptsize $b$}
\psfrag{a1}{\scriptsize $a'$}
\psfrag{b1}{\scriptsize $b'$}
\psfrag{a2}{\scriptsize $a'$}
\psfrag{b2}{\scriptsize $b'$}
\psfrag{aa}{\small (a)}
\psfrag{bb}{\small (b)}
\psfrag{cc}{\small (c)}
\psfrag{i}{\scriptsize $a_i=b_i$}
\epsfig{file=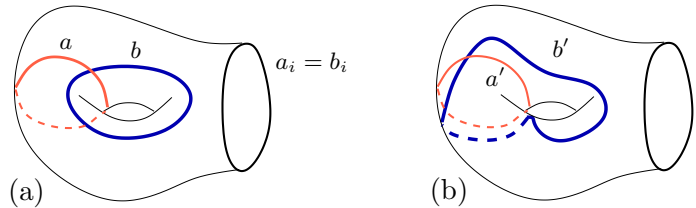,width=0.59\linewidth}
\caption{Handle-twists:
(a) Double self-folded pair of pants; (b) The same pair of pants after a handle-twist $T_a(b)$} 
\label{d-self-pant}
\end{center}
\end{figure}

\begin{definition}[{{\it Handle-twists}}]
\label{Dehn}
Given a double pants decomposition $DP=(P_a,P_b)$
we define an additional transformation which may be performed if $P_1$ and $P_2$ contain the same  curve $a_i=b_i$  
separating the same handle $\h$, 
see Fig.~\ref{d-self-pant}(b). Let $a\in \h$ and $b\in \h$ be the only curves from $P_a$ and $P_b$ respectively.
Then a {\it handle-twist} $T_a(b)$ (respectively, $T_b(a)$) is a Dehn twist along $a$  (respectively, $b$) in any of two directions
(see Fig.~\ref{d-self-pant}(b)).  

\end{definition}

Notice that both flips and handle-twists are reversible transformations, so that flips and handle-twists generate a groupoid acting on
the set of double pants decompositions.

\begin{definition}[{{\it Flip-twist groupoid}}]
A {\it flip-twist groupoid} $FT$ is a groupoid generated by flips and twists. 

\end{definition}

\begin{definition}[{{\it Double curve}}]
A curve $c\in (P_a,P_b)$ is {\it double} if $c\in P_a\cap P_b$. 

\end{definition}

\begin{definition}[{{\it Standard decomposition, principle curves}}]
A double pants decomposition $(P_a,P_b)$ is {\it standard} if 
there exist $g$ double curves $c_1,\dots,c_g \in (P_a,P_b)$ such that
 $c_i$ cut out of $S$ a handle $\h_i$. 

\end{definition}




\begin{figure}[!h]
\begin{center}
\psfrag{a}{$P_a$}
\psfrag{b}{$P_b$}
\epsfig{file=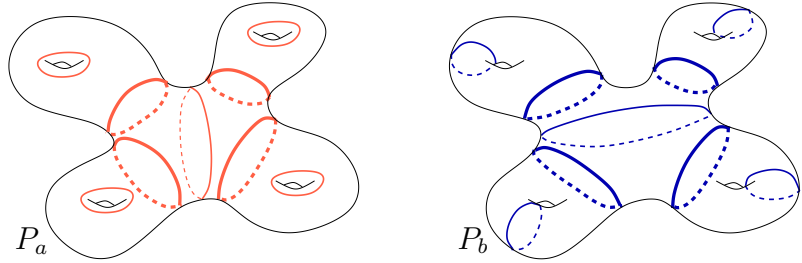,width=0.68\linewidth}
\caption{A standard double pants decomposition $(P_a,P_b)$.  } 
\label{standard double_pant}
\end{center}
\end{figure}

%

\begin{definition}[{{\it Admissible decomposition}}]
\label{admissible}
A double pants decomposition $(P_a,P_b)$ is {\it admissible}  if it is possible to transform $(P_a,P_b)$ to a standard pants 
decomposition by a sequence of flips.


\end{definition} 

%


%

The following theorem is the main result of~\cite{FN}.

\begin{theorem}[\cite{FN}]
\label{unlabeled} 
A flip-twist groupoid acts transitively on admissible double pants decompositions of $S=S_{g,n}$ (for any $(g,n)$ such that $2g+n>2$).

\end{theorem}

%

\subsection{Labeled double pants decompositions}
A pants decomposition\! $P\!=\!\langle c_1,\!\dots\!,c_m\rangle$ is {\it  labeled} if each curve  $c_i\in P$ is labeled by a number 
$x_i\in \{1,\dots,m\}$, $x_i\ne x_j$ for $i\ne j$, in other words we assign to the curves of $P$ distinct numbers $\{1,\dots, m\}$.

Similarly, 
a double pants decomposition $DP=(P_a,P_b)=\langle c_1,\dots,c_m;c_{m+1},\dots,c_{2m}\rangle$ is {\it  labeled} if each curve  
$c_i\in P_a$ is labeled by a number $x_i\in \{1,\dots,2m\}$.

\begin{figure}[!h]
\begin{center}
\psfrag{1}{\scriptsize $1$}
\psfrag{2}{\scriptsize $2$}
\psfrag{3}{\scriptsize $3$}
\psfrag{a}{\small (a)}
\psfrag{b}{\small (b)}
\epsfig{file=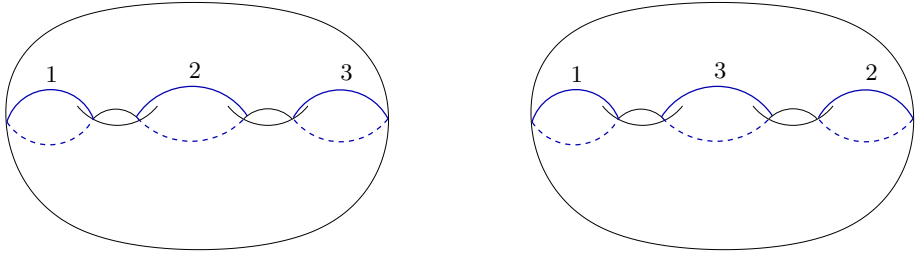,width=0.78\linewidth}
\caption{These two labeled decompositions are different.} 
\label{labeled}
\end{center}
\end{figure}

In unlabeled version of double pants decomposition we consider $\langle c_1,\dots,c_m\rangle$ as a set, so we do not distinguish between
two pants decompositions shown in Fig.~\ref{labeled}. In labeled version this two decompositions are considered as different
(they differ by their {\it labelings}).

The flips and handle-twists preserve the labelings, i.e. the new curve curries the same number as the deleted one had.
In case of a double curve $c=c_i=c_j\in DP$ we assign two labels $x_i$ and $x_j$ to the same curve. 
Flipping the double curve we can not differ between 
the labels, so for any given topological flip $f$ of a double curve we define two labeled flips $f_1$ and $f_2$:
both result in the same set of curves on the surface as $f$ does, but the curve $c$ is labeled by $x_i$ after $f_1$   
and by $x_j$ after $f_2$.

In~\cite{FN} we have proved transitivity of flip-twist groupoid on admissible double pants decompositions. It is natural to ask if
this groupoid act transitively on labeled admissible double pants decompositions.

\subsection{Strict-labeled double pants decompositions}

In a labeled double pants decomposition the labels of a double curve are fully equal in rights. 
It would be convenient for the proof of our theorem to define also {\it strict-labeled} double pants decompositions,
where each label of a double curve remembers to which of two pants decompositions it belongs (so that flips and twists
do not intermix labels of $P_a$ with labels of $P_b$.

A double pants decomposition $DP=(P_a,P_b)=\langle c_1,\dots,c_m;c_{m+1},\dots,c_{2m}\rangle$ is {\it  strict-labeled} if each curve  
$c_i\in P_i$ is labeled by a number $x_i\in \{1,\dots,m\}$ and each curve $c_i\in P_j$  is labeled by a number 
$x_i\in \{m+1,\dots,2m\}$, where either $P_i=P_a$, $P_j=P_b$ or  $P_i=P_b$, $P_j=P_a$. 
In other words, instead of the whole permutation group $S_{2m}$ the set of labels is permuted only by $S_m\times \Z_2$.

The flips and handle-twists preserve the labeling, i.e. the new curve curries the same number as the deleted one had.
In case of a double curve $c=c_i=c_j\in DP$ we  do not mix $c_i$ with $c_j$ since one of these curves belong to $P_a$ and another belong
to $P_b$. This mean that if $c_i\in P_a$ and we make flip $f$ of $c\in P_a$, then the new curve $f(c)\in P_a$  
has the same label as $c_i$, not as $c_j$. So the labels assigned to one component always stay together.

We will first work with strict-labeled double pants decompositions and then in Section~\ref{sec weak} extend the results
to the labeled double pants decompositions.

\section{Transitivity for labeled  pants decompositions}
\label{sec ordinary}

\begin{figure}[!h]
\begin{center}
\psfrag{1}{\scriptsize $c$}
\psfrag{a}{\scriptsize $c_1$}
\psfrag{a1}{\scriptsize $c_1'$}
\epsfig{file=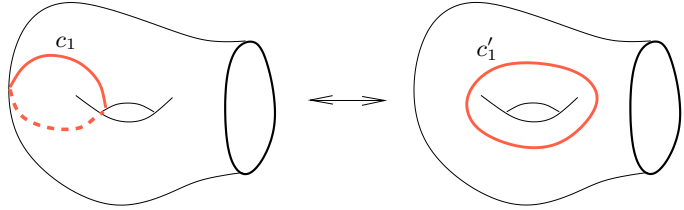,width=0.58\linewidth}
\caption{$S$-move: if some handle $\h$ is separated by by a curve $c\in P$ then a curve $c_1\in P$ contained in $\h$ may be exchanged by
any curve $c_1'$ such that $|c_1\cap c_1'|=1$.} 
\label{s_move}
\end{center}
\end{figure}

Let $P$ be a pants decomposition of a surface $S=S_{g,n}$
It is shown in~\cite{HT} and~\cite{H1}   that $P$ may be transformed to any other pants decomposition of $S$
via a sequence of flips and $S$-moves, where an $S$-move is defined as in Fig.~\ref{s_move}  
In this section we show that the groupoid generated by flips and $S$-moves acts transitively on labeled pants decompositions.

First we will prove transitivity for the case of sufficiently large surfaces, namely for $S_{g,n}$ satisfying $2g+n>4$,
or in other words, to the surfaces whose pants decomposition contain at least three pairs of pants.

\begin{lemma}
\label{neighbours}
Let $P$ be a pants decomposition of a surface $S_{g,n}$, $2g+n>4$. 
If $c_1,c_2\in P$ are two curves in the boundary of the same 
pair of pants then the label of $c_1$ may be swapped with the label of $c_2$  by a sequence of flips.

\end{lemma}

\begin{proof}
Let $p$ be a pair of pants containing both $c_1$ and $c_2$ as boundary components. Let $p_1$ and $p_2$ be the adjacent pairs of pants
($c_i=p\cap p_i$). There are two possibilities: either $p_1$ and $p_2$ are two distinct pairs of pants or they coincide.

If $p_1\ne p_2$ then  the labels on $c_1$ and $c_2$ are swapped by a sequence of 5 flips shown in Fig.~\ref{pentagon}
(usually called ``pentagon relation''). Notice that some of five boundary components of the ``pentagon'' may be identified,
but this does not affect the procedure.

\begin{figure}[!h]
\begin{center}
\psfrag{1}{\scriptsize $1$}
\psfrag{2}{\scriptsize $2$}
\psfrag{3}{\scriptsize $3$}
\psfrag{a}{\small (a)}
\psfrag{b}{\small (b)}
\epsfig{file=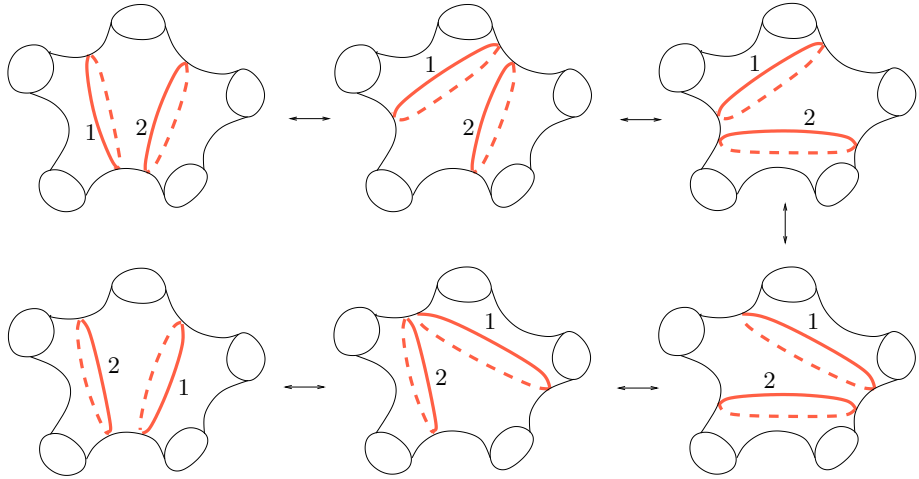,width=0.78\linewidth}
\caption{``Pentagon relation'' exchanges the labels} 
\label{pentagon}
\end{center}
\end{figure}

Suppose that $p_1=p_2$. Since $S$ contains at least 3 pairs of pants, there exists a pair of pants $p_3$ such that $p\cup p_1\cup p_3$
looks as shown in Fig.~\ref{3pants}, left  (up to possible interchange of $p$ and $p_1$ and possible identification of some boundary components). 
Then after one flip we obtain a configuration on 
 Fig.~\ref{3pants}, right, which suits to the case $p_1\ne p_2$ considered above.

\begin{figure}[!h]
\begin{center}
\psfrag{1}{\scriptsize $1$}
\psfrag{2}{\scriptsize $2$}
\psfrag{p}{\scriptsize $p$}
\psfrag{p1}{\scriptsize $p_1=p_2$}
\psfrag{p_1}{\scriptsize $p_1$}
\psfrag{p_2}{\scriptsize $p_2$}
\psfrag{p3}{\scriptsize $p_3$}
\psfrag{a}{\small (a)}
\psfrag{b}{\small (b)}
\epsfig{file=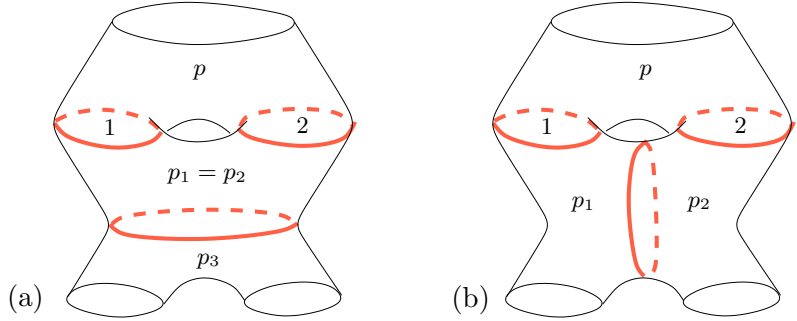,width=0.68\linewidth}
\caption{Reduction of the case  $p_1=p_2$ to the case  $p_1\ne p_2$} 
\label{3pants}
\end{center}
\end{figure}

\end{proof}

\begin{cor}
\label{large}
Let $P$ be a labeled pants decomposition of $S=S_{g,n}$ where
 $2g+n>4$. Then flips act transitively on labeling of $P$.

\end{cor}

\begin{proof}
The statement follows immediately from Lemma~\ref{neighbours} and the fact that $S$ is connected.

\end{proof}

\begin{figure}[!h]
\begin{center}
\psfrag{3}{\scriptsize $1$}
\psfrag{4}{\scriptsize $2$}
\psfrag{S3}{\scriptsize $S_1$}
\psfrag{F4}{\scriptsize $F_2$}
\epsfig{file=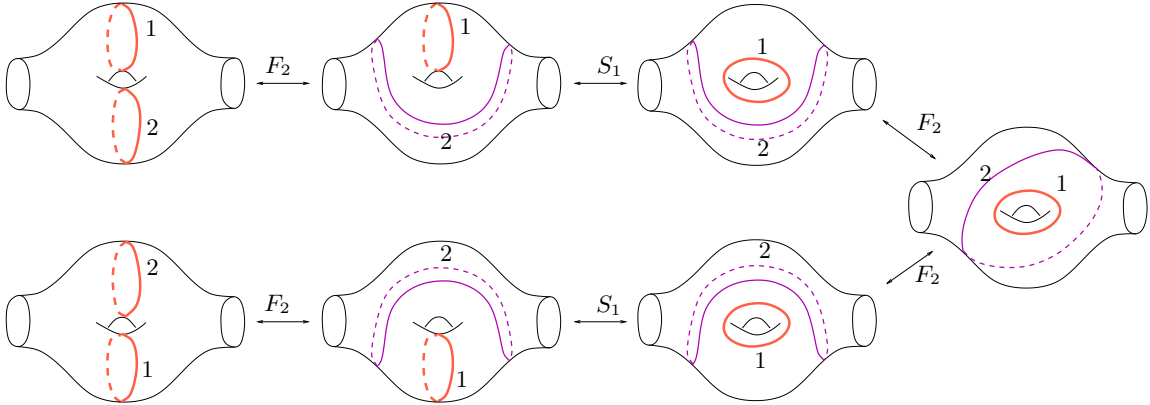,width=0.98\linewidth}
\caption{Exchange of labels in $S_{1,2}$} 
\label{s-1_2}
\end{center}
\end{figure}

It is clear that Corollary~\ref{large} together with transitivity of flips and $S$-moves on unlabeled pants decompositions
imply transitivity on  the labeled pants decompositions. 
So, we are left to consider only finitely many surfaces $S_{g,n}$ satisfying the inequality $2g+n\le 4$,
i.e. $S_{0,3}$, $S_{0,4}$, $S_{1,1}$, $S_{1,2}$ and $S_{2,0}$. Notice, that pants decompositions of surfaces
 $S_{0,3}$, $S_{0,4}$ and $S_{1,1}$ contain at most one curve, so there is nothing to prove in these cases.
A pants decomposition of $S_{1,2}$ contains two curves whose labels could be swapped as in Fig.~\ref{s-1_2}.
So, the only question left concerns $S_{2,0}$. 
A pants decomposition of $S_{2,0}$ contains 3 curves. After at most one flip we may assume that all three curves are homologically 
non-trivial. Then we cut $S_{2,0}$ along one of the curves and use the procedure described above for $S_{1,2}$ to swap the labels of 
the other two curves. Thus, we have all transpositions of the labels on $P$, and hence, all permutations.

We summarize the results of this section in the following theorem:

\begin{theorem}
For any surface $S_{g,n}$, $2g+n>2$ flips and $S$-moves act transitively on labeled pants decompositions of $S$.

\end{theorem}

\section{Transitivity for strict-labeled double pants decompositions}
\label{sec strict}

Our next step is to prove transitivity of flips-twist groupoid on the strict-labeled double pants decompositions.

\begin{theorem}
\label{large double}
Let $DP$ and $DP'$ be two strict-labeled admissible double pants decompositions of  $S=S_{g,n}$ where $2g+n>4$. 
Then there exists a sequence of flips and handle-twists transforming $DP$ to $DP'$.

\end{theorem}

\begin{proof}
By Theorem~\ref{unlabeled} there exists a sequence $\psi$ of flips and handle-twists which takes $DP=(P_a,P_b)$ to $DP'=(P_a',P_b')$ 
as an unlabeled decomposition. There are two possibilities: either $\psi$ takes $P_a$ to $P_a'$ and $P_b$ to $P_b'$
or $\psi$ changes the components $P_a$ and $P_b$.

Suppose that $\psi(P_a)=P_a'$, $\psi(P_b)=P_b'$. Then we apply Corollary~\ref{large} to take labeling of $\psi(P_a)$ and $\psi(P_b)$
to labeling of $P_a'$ and $P_b'$ respectively.

Suppose that  $\psi(P_a)=P_b'$, $\psi(P_b)=P_a'$. Consider a sequence of flips and handle-twists $\varphi$ which takes $(P_a',P_b')$ 
to a standard pants decomposition (it does exists by  Theorem~\ref{unlabeled}). In a standard double pants decomposition 
we may swap the
curve of $P_a$ with the curve of $P_b$  in each handle separately: see Fig.~\ref{handle S-move}.

\begin{figure}[!h]
\begin{center}
\psfrag{a}{\scriptsize $i$}
\psfrag{b}{\scriptsize $j$}
\epsfig{file=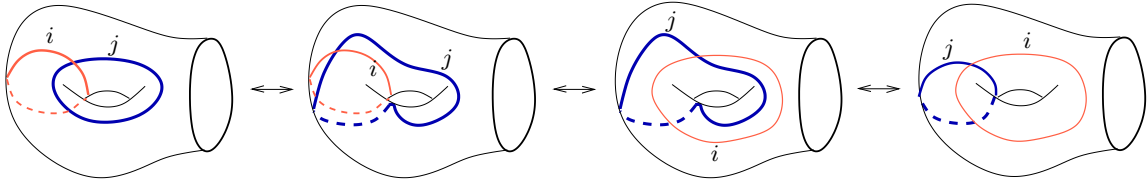,width=0.98\linewidth}
\caption{Exchange of labels in a handle via 3 handle-twists} 
\label{handle S-move}
\end{center}
\end{figure}

\end{proof}


In view of Theorem~\ref{large} we are left to consider only finitely many surfaces $S_{g,n}$ satisfying the inequality $2g+n\le 4$,
i.e. $S_{0,3}$, $S_{0,4}$, $S_{1,1}$, $S_{1,2}$ and $S_{2,0}$. Below, we consider these five surfaces one by one.

\subsection{Surface $S_{0,3}$.}
\label{S_{0,3}}
A pants decomposition of this surface is trivial (contains no curve), so the transitivity of flip-twist groupoid on 
strict-labeled pants 
decompositions follows trivially.

\subsection{Surface $S_{0,4}$.}
A pants decomposition of the sphere with 4 holes consists of one curve, a double pants decomposition is an arbitrary  pair of curves
$\langle a,b \rangle$.
In particular, after one flip $f$  we may assume that the $a$ and $ f(b)$ coincide, and the inverse of this flip applied to $a$
gives $\langle f^{-1}a,f(b)\rangle =\langle b,a\rangle$. This shows that applying flips we may swap the label of $a$ with the label of $b$.

\subsection{Surface $S_{1,1}$.}
A pants decomposition of the sphere with 4 holes consists of one curve, a double pants decomposition is a   pair of curves
$\langle a,b \rangle$, such that $|(a,b)|=1$. The labels of these curves may be swapped by a sequence of three handle-twists, as in 
 Fig.~\ref{handle S-move}.

\subsection{Surface $S_{1,2}$.}
\begin{figure}[!h]
\begin{center}
\psfrag{1}{\scriptsize $1$}
\psfrag{2}{\scriptsize $2$}
\psfrag{3}{\scriptsize $3$}
\psfrag{4}{\scriptsize $4$}
\psfrag{2,3}{\scriptsize $2,3$}
\psfrag{2,4}{\scriptsize $2,4$}
\psfrag{f2}{\scriptsize $F_2$}
\psfrag{f3}{\scriptsize $F_3$}
\psfrag{f4}{\scriptsize $F_4$}
\psfrag{T3}{\scriptsize $T_3$}
\psfrag{T4}{\scriptsize $T_4$}
\epsfig{file=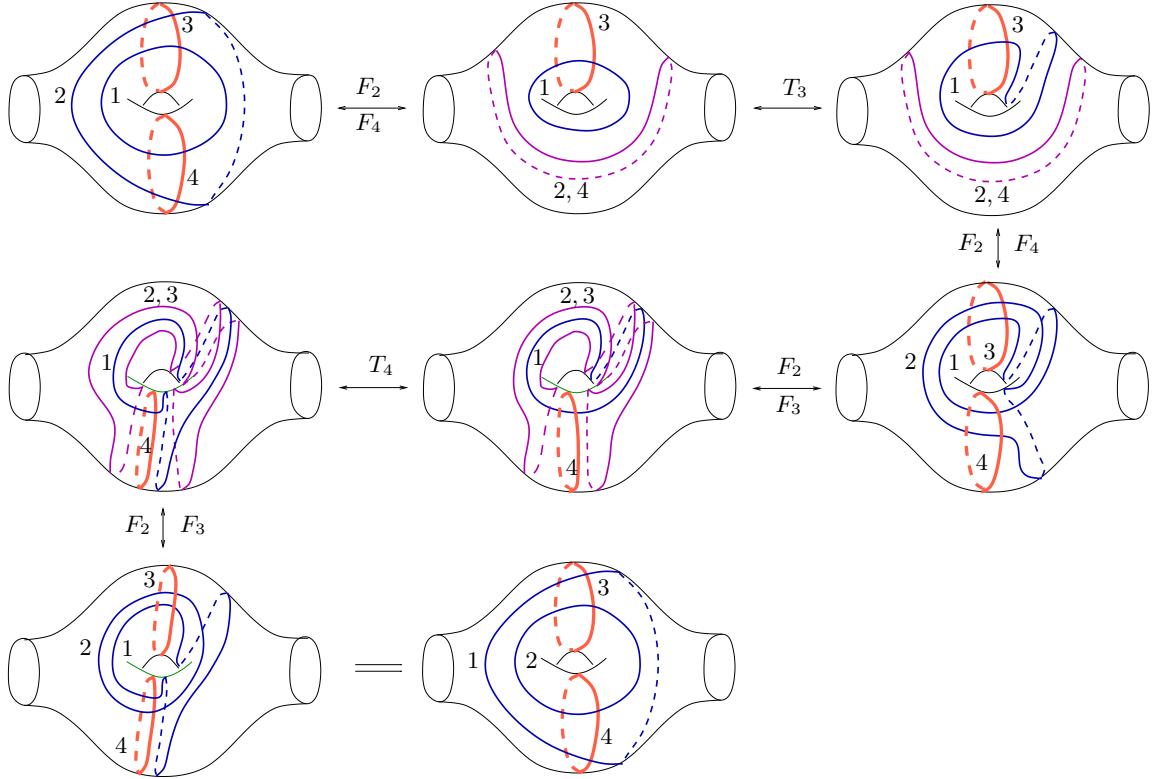,width=0.98\linewidth}
\caption{Changing the labels of $P_a$ and preserving the labels of $P_b$ ($F_i$ stays for a flip of $i$-th curve, $T_i$ stays for a 
handle-twist along $i$-th curve)
} 
\label{S_1_2}
\end{center}
\end{figure}

A pants decomposition of $S_{1,2}$ consists of two curves, a double pant decomposition consists of two pairs of curves.
The labels of the two components $P_a$ and $P_b$ may be swapped in the usual way: in a standard double pant decomposition
we change the labels via three handle-twists, as in Fig.~\ref{handle S-move}. So, the only thing to check
is that there exists a sequence of flips and handle-twists which exchanges the labels of two curves of $P_a$ and preserves the 
labels of $P_b$. This sequence is shown in Fig.~\ref{S_1_2}.

\subsection{Surface $S_{2,0}$.}
\begin{figure}[!h]
\begin{center}
\psfrag{a1}{\scriptsize $a_1$}
\psfrag{a2}{\scriptsize $a_2$}
\psfrag{a3}{\scriptsize $a_3$}
\psfrag{b1}{\scriptsize $b_1$}
\psfrag{b2}{\scriptsize $b_2$}
\psfrag{b3}{\scriptsize $b_3$}
\epsfig{file=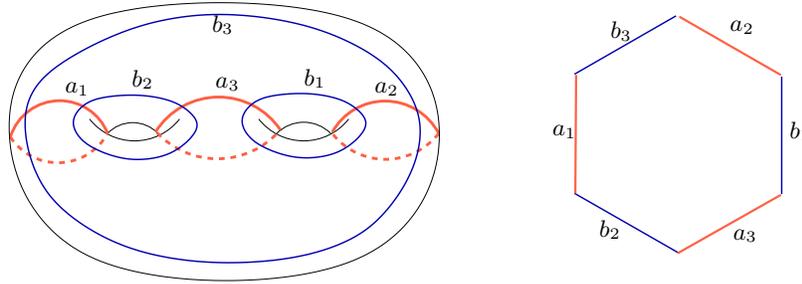,width=0.68\linewidth}
\caption{A hexagonal double pants decomposition} 
\label{ideal}
\end{center}
\end{figure}

\label{S_{2,0}}
This is the only surface where flip-twist groupoid do not act transitively on strict-labeled admissible double pants decompositions.
To investigate this action we consider the double pants decompositions of the combinatorial type shown in Fig.~\ref{ideal} 
(we call this type of double pants decomposition {\it hexagonal}).
This type of double  pants decomposition may be characterized by the following property: 
we have $|a_i\cap b_{i+1}|=1$,  $|b_i\cap a_{i+1}|=1$, $|a_i\cap b_i)|=0$, where indexes are considered modulo 3. 
In other words, the curves $[a_1,b_3 ,a_2 ,b_1,a_3 ,b_2]$ compose ``hexagon'', where the neighbours do intersect and other sides 
do not. More precisely, the curves of the hexagonal decomposition decompose the surface $S_{2,0}$ into 4 hexagons.

Our aim now is to label the curves $[1,4,2,5,3,6]$ and check which permutations of these labels may be performed via flips and
handle-twists.

Consider any of the four hexagons on $S_{2,0}$ and read the labels on its sides in a clockwise direction. The obtained sequence
should be considered modulo cyclic shifts  of all labels and modulo reversing of the order.
We obtain a {\it cyclic order} of the labeling.
Notice that the cyclic order does not depend on the choice of one of the four hexagons 
(the choice of any of the two adjacent hexagons reverse the order of the labels in the sequence, the choice of the hexagon opposite 
to the initial one does not affects the order).

Whenever we need to compare two labeling of the same hexagonal double pants decomposition we always refer to the same hexagon on the 
surface, so that the notions of rotation of the labels and reversing of the order of the labels make sense.
Comparing the cyclic orders of two distinct hexagonal decompositions we think of the cyclic order modulo reversing of the order
and cyclic shifts of the labels.

\begin{definition}[{{\it Hexagonal twist}}]
Let  $[a_1,b_3 ,a_2 ,b_1,a_3 ,b_2]$ be a hexagonal double pants decomposition. A {\it hexagonal twist}  $T_{a_i}$ (or $T_{b_i}$) is a
Dehn twist along $a_i$ (respectively, along $b_i$), $i=1,2,3$.

\end{definition}

\begin{lemma}
\label{hexagonal twist}
\begin{itemize}
\item[1)] Any hexagonal twist is a composition of flips and handle-twists;
\item[2)] Any handle-twist is a composition of flips and hexagonal twists;
\item[3)] A hexagonal twist preserves the cyclic order in a hexagonal double pants decomposition.

\end{itemize}
\end{lemma}

\begin{figure}[!h]
\begin{center}
\psfrag{a1}{\scriptsize 1}
\psfrag{a2}{\scriptsize 2}
\psfrag{a3}{\scriptsize 3}
\psfrag{b1}{\scriptsize 4}
\psfrag{b2}{\scriptsize 5}
\psfrag{b3}{\scriptsize 6}
\psfrag{a3,b3}{\scriptsize 3,6}
\psfrag{hex}{\scriptsize Hexagonal twist}
\psfrag{han}{\scriptsize Handle-twist}
\psfrag{1}{\scriptsize flip of $a_3$,}
\psfrag{2}{\scriptsize flip of $b_3$}
\epsfig{file=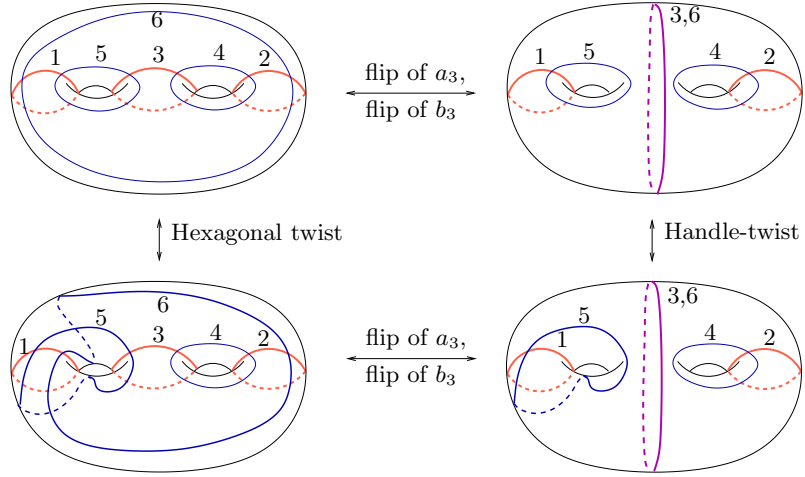,width=0.68\linewidth}
\caption{Hexagonal twist as a composition of flips and handle-twists} 
\label{hex is handle}
\end{center}
\end{figure}

\begin{proof}
Parts 1) and 2) follow from the commutative diagram shown in Fig~\ref{hex is handle},
part 3) is evident.


\end{proof}

\begin{lemma}
\label{l rotation}
There exists a sequence of flips and handle-twists which takes the set of curves  $[a_1,b_3 ,a_2 ,b_1,a_3 ,b_2]$ to itself and the labels
 $[1,4,2,5,3,6]$ to  $[6,1,4,2,5,3]$.

\end{lemma}

\begin{proof}
Consider the system  $[a_1,b_3 ,a_2 ,b_1,a_3 ,b_2]$ and apply five hexagonal twists in a row, namely  $T_2$, $T_4$, $T_3$, 
$T_5$, $T_1$, where $T_i$ is a twist in the curve labeled by $i$.
Then we return to the same set of curves, but the labels are shifted, see Fig.~\ref{rotation}.

\begin{figure}[!h]
\begin{center}
\psfrag{a1}{\scriptsize $1$}
\psfrag{a2}{\scriptsize $2$}
\psfrag{a3}{\scriptsize $3$}
\psfrag{b1}{\scriptsize $4$}
\psfrag{b2}{\scriptsize $5$}
\psfrag{b3}{\scriptsize $6$}
\psfrag{Ta1}{\scriptsize $T_1$}
\psfrag{Ta2}{\scriptsize $T_2$}
\psfrag{Ta3}{\scriptsize $T_3$}
\psfrag{Tb1}{\scriptsize $T_4$}
\psfrag{Tb2}{\scriptsize $T_5$}
\psfrag{Tb3}{\scriptsize $T_6$}
\epsfig{file=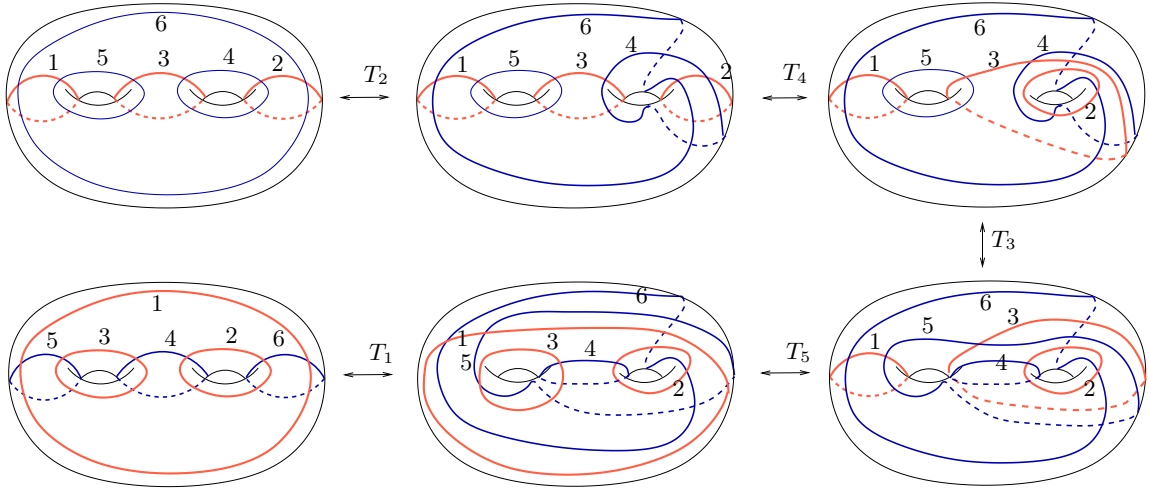,width=0.98\linewidth}
\caption{Rotation of the hexagon realized by 5 hexagonal twists.} 
\label{rotation}
\end{center}
\end{figure}

\end{proof}

\begin{lemma}
\label{l reflection}
There exists a sequence of flips and handle-twists which takes the set of curves  $[a_1,b_3 ,a_2 ,b_1,a_3 ,b_2]$ to itself and the labels
 $[1,4,2,5,3,6]$ to  $[6,3,5,2,4,1]$.

\end{lemma}

\begin{figure}[!h]
\begin{center}
\psfrag{f1}{\scriptsize $F_3$}
\psfrag{f2}{\scriptsize $F_6$}
\psfrag{a1}{\scriptsize $1$}
\psfrag{a2}{\scriptsize $2$}
\psfrag{a3}{\scriptsize $3$}
\psfrag{b1}{\scriptsize $4$}
\psfrag{b2}{\scriptsize $5$}
\psfrag{b3}{\scriptsize $6$}
\psfrag{Ta1}{\scriptsize $T_1$}
\psfrag{Ta2}{\scriptsize $T_2$}
\psfrag{Ta3}{\scriptsize $T_3$}
\psfrag{Tb1}{\scriptsize $T_4$}
\psfrag{Tb2}{\scriptsize $T_5$}
\psfrag{Tb3}{\scriptsize $T_6$}
\epsfig{file=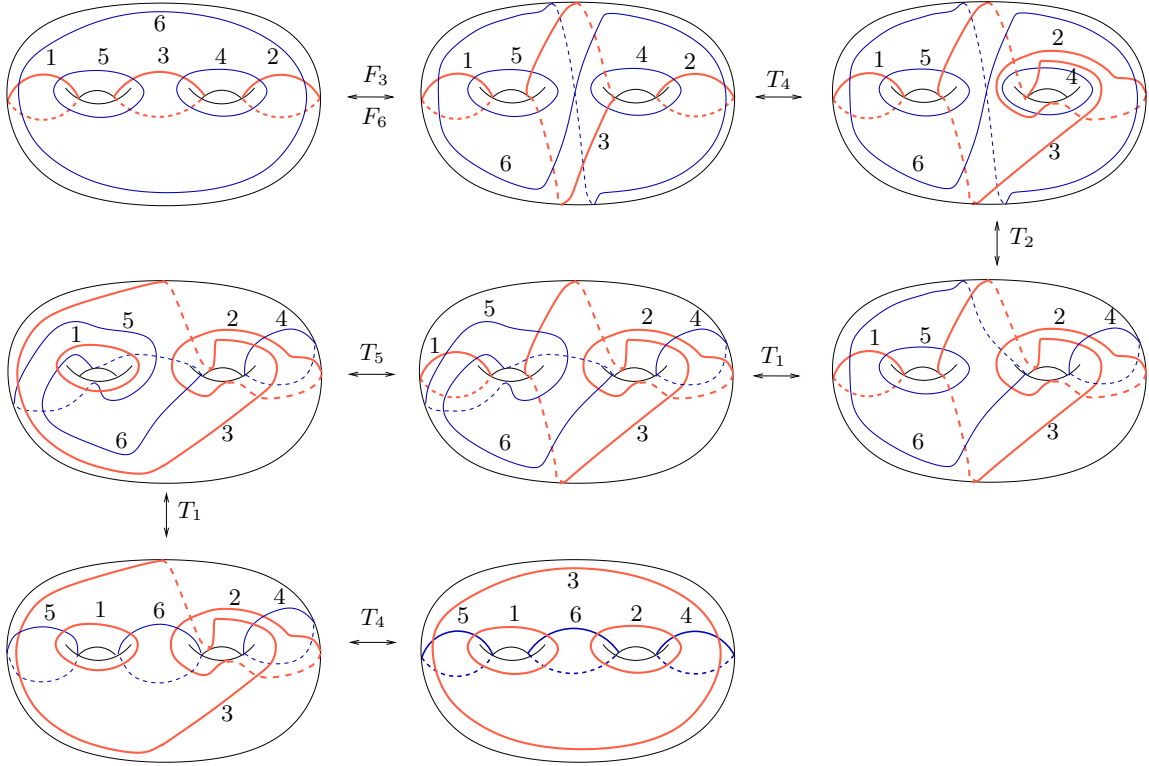,width=0.98\linewidth}
\caption{Reflection of the hexagon  $[a_1,b_3 ,a_2 ,b_1,a_3 ,b_2]$.} 
\label{reflection}
\end{center}
\end{figure}

\begin{proof}
First, we make a flip in 3 and 6, and then apply  hexagonal twists in 4, 2, 1, 5, 1, 4,
see Fig.~\ref{reflection}

\end{proof}

Lemmas~\ref{l rotation} and~\ref{l reflection} show that show that each element of the dihedral group $D_6$ acting on the labels 
of hexagon
$[a_1,b_3 ,a_2 ,b_1,a_3 ,b_2]$ may be realized as a sequence flips and handle-twists (Lemmas~\ref{l rotation} and~\ref{l reflection}
 represent a rotation and reflection of the hexagon respectively).

\begin{remark}
We generate $D_6$ be rotation of order 6 and one reflection. Another possibility is to generate $D_6$ by two reflections, but then we 
need to find a presentation (via flips and twists) of a reflection saving two opposite sides of the hexagon. 

\end{remark}

In  Lemma~\ref{dihedral} we will show that there are no other permutations of labels of the hexagon that can be realized by flips and  
handle-twists. For the proof we will consider $\Z_2$ homology classes.
Denote $\bar h(c)\in H_1(S,\Z_2)$ a $\Z_2$ homology class of the curve $c$.
Notice that a homology class $h(c)\in H_1S,\Z)$ is defined only up to a change of sign (depending on the orientation of the curve $c$),
however, the class $\bar h(c)\in H_1(\Z_2,S)$ is well defined for a non-oriented curve.

\begin{lemma}
\label{Z_2 homology}
Let $DP=[a_1,b_3 ,a_2 ,b_1,a_3 ,b_2]$ be a hexagonal double pants decomposition of $S_{2,0}$.
Let $\psi$ be a sequence of flips of $P_a=\langle a_1,a_2,a_3\rangle$ such that $\psi(P_a)$ contains no homo logically trivial curves.
Then $\bar h(\psi(a_i)=\bar h(a_i)$, $i=1,2,3$.  

\end{lemma}

\begin{proof}
Label the curves $a_1,a_2,a_3$ of $P_a$ by numbers $1,2,3$ respectively and
consider the sequence $\psi$ as a composition of subsequences $\psi_1\circ\dots\circ\psi_k$, where   $\psi_i$ is a composition of
flips of the curve with the same label, while $\psi_i$ and $\psi_{i+1}$ are compositions of flips of curves with different labels.
It is easy to see that after applying any subsequence $\psi_1\circ\dots\circ\psi_j$, $0\le j\le k$ the pants decomposition $P_a$ 
turns in a decomposition without homologically trivial curves (otherwise $\psi_j$ and  $\psi_{j+1}$ can not flip the curve with 
distinct labels). So, it is sufficient to prove the statement of the lemma for one subsequence $\psi_i$.

Now, suppose that all flips in $\psi_1$ change the curve labeled $1$. Consider a pair of pants $p$ in $\psi_1(P_a)$ containing the curve
$\psi_1(a_1)$. Then the boundary of $p$ consists of curves $\psi_1(a_1)$, $a_2$ and $a_3$ 
(notice that $\partial p\ne \psi_1(a_1)\cup a_i$, $i=2$ or $3$ since $\psi_1(a_1)$ is homologically non-trivial).
This implies that $h(\psi_1(a_1))=\pm (h(a_2) \pm h(a_3))$
(the choices of signs depend on the orientations of curves). So, for $\Z_2$-homology classes we always have 
 $$\bar h(\psi_1(a_1))=\bar h(a_2) + \bar h(a_3))=\bar h(a_1),$$
and the statement for $\psi_1$ is proved. Applying this $k$ times we obtain the lemma.

\end{proof}

%
%


\begin{lemma}
\label{dihedral}
Let $\varphi$ be a sequence of flips and handle-twists transforming the hexagonal set of curves $[a_1,b_3 ,a_2 ,b_1,a_3 ,b_2]$ 
to itself and 
permuting the labels $[1,4,2,5,3,6]$ of these curves. Then the permutation coincides with some permutation obtained by an action of
dihedral group $D_6$ on the hexagon  $[a_1,b_3 ,a_2 ,b_1,a_3 ,b_2]$.

\end{lemma}

\begin{proof}
Consider a sequence $\varphi$ of flips and handle-twists transforming the set of curves  $\langle a_1,b_3 ,a_2 ,b_1,a_3 ,b_2\rangle$ to itself.
By Lemma~\ref{hexagonal twist} any handle-twist is a composition of flips and hexagonal twists. So, $\varphi$ is a composition of flips 
and hexagonal twists. 
By Lemma~\ref{hexagonal twist} hexagonal twist does not changes cyclic order.  
We will show  that if a sequence of flips takes a hexagonal double pants decomposition to a hexagonal one, then it 
either preserves the cyclic order or changes it to the opposite. Then the  statement of the lemma follows.

So, we are left to show that if a sequence of flips $\varphi$ takes a hexagonal double pants decomposition to a hexagonal one, 
then $\varphi$ either does not change the cyclic order or .
By Lemma~\ref{Z_2 homology}, $\varphi$ preserves $\Z_2-$homology classes of all curves of $P_a$ and $P_b$. This implies that 
$\varphi(a_i)$ intersects $\varphi(b_j)$ if and only if $a_i$ intersects $b_j$. So, the cyclic order of curves 
$[a_1,b_3 ,a_2 ,b_1,a_3 ,b_2]$ is either preserved by $\varphi$ ot changed to the opposite.
 
\end{proof}

\begin{cor}
The action of flip and twist groupoid on strict-labeled admissible double pants decompositions of $S_{2,0}$ has  6 orbits.

\end{cor}

\begin{proof}
Consider a hexagonal double pants decomposition.
Notice that the labels from the set $\{1,2,3\}$ alternate with the labels from the set $\{3,4,5\}$, since curves from the same pants 
decomposition do not intersect each other.  

Now, the action of the dihedral group $D_6$ is sufficient to put the labels $1,2,3$ of $P_a$ into required position.
There are 6 possibilities left to put the labels $4,5,6$. Distinct possibilities result in distinct cyclic orders, so the obtained 
labeling are not equivalent under the action of flip and twist groupoid.

\end{proof}

The results of Section~\ref{S_{0,3}}--Section~\ref{S_{2,0}} may be summarized in the following theorem:

\begin{theorem}
\label{trans}
The flip-twist groupoid acts transitively on strict-labeled admissible double pants decompositions of $S_{g,n}$ unless $(g,n)=(2,0)$.
The action of flip and twist groupoid on strict-labeled admissible double pants decompositions of $S_{2,0}$ has  6 orbits.

\end{theorem}

\section{Transitivity on labeled double pants decompositions}
\label{sec weak}
Now, we will extend Theorem~\ref{trans} from the class of strict-labeled decompositions to the class of labeled decompositions.

Instead of thinking that two labels $x_1$ and $x_2$ of a double curve $c$ do not know  to which of $P_a$ and $P_b$ they do belong,
we will add one more transformation which we will call {\it switch}: a switch in $c$ just changes the places of $x_i$ and $x_j$.


\begin{theorem}
\label{weak}
The flip-twist groupoid acts transitively on labeled admissible double pants decompositions of $S_{g,n}$ unless $(g,n)=(2,0)$.
The action of flip and twist groupoid on labeled admissible double pants decompositions of $S_{2,0}$ has  $15$ orbits.

\end{theorem} 

\begin{proof}
First, consider all surfaces $S_{g,n}$, $(g,n)\ne (2,0)$. 
By Theorem~\ref{trans} flip-twist groupoid act transitively on strict-labeled double pants decompositions of these surfaces. So, it is 
sufficient to check that applying flips, handle-twists and switches in double curves we may arrive in a double pants decomposition 
$DP=(P_a,P_b)$ where the labels of $P_a$ contain any given set of $m$ numbers from  $\{1,2,\dots,2m\}$.
Clearly, it is sufficient to check that we may change the places of two labels, assigned to curves in $P_a$ and $P_b$ respectively.  
More precisely, if $c_i\in P_a$, $c_j\in P_b$, then we may change the labels $i$ and $j$ of the curves $c_i$ and $c_j$ in the 
following way:
\begin{itemize}
\item[1)] find an admissible decomposition $DP'$ with a double curve $c$ (we can get to $DP'$ from $DP$ by flips and handle-twists
in view of Theorem~\ref{unlabeled});
\item[2)] use transitivity of flip-twist groupoid on labeled double pants decompositions to get to the decomposition $DP'$ with labels
$i$ and $j$ on the double curve $c$;
\item[3)] switch $i$ and $j$;
\item[4)] return to the decomposition $DP$ with $c_i$ and $c_j$ labeled $j$ and $i$ respectively.  

\end{itemize}

This implies that if   $(g,n)\ne (2,0)$ then the flip-twist groupoid acts transitively on labeled admissible
double pants decompositions.

\medskip

Now consider the case   $(g,n)=(2,0)$.  
The switch of the labels on the double curve in terms of hexagonal decomposition describes as exchange of the labels of two opposite
sides of the hexagon. This implies that under the action of flip-twist groupoid the labels of opposite sides always remain the opposite. 
There are 15 possibilities to split $6$ labels into 3 pairs (the label 1 may be in pair with each of 5 other labels, for each of these 
5 possibilities the smallest of the remaining labels may be paired with any of 3 other labels).

Choose one of the possible pairings, for example, $\langle(1,4),(2,5),(3,6)\rangle$. Using switches we make the labels $1,2,3$ mutually
non-adjacent in the hexagon. Using rotations and reflections as in Lemmas~\ref{l rotation} and~\ref{l reflection}
we may put the labels $1,2,3$ in any three mutually non-adjacent positions. The pairing then determines the positions of all remaining 
labels.  

\end{proof}


\begin{thebibliography}{5}	
\bibitem{C} Cameron McA Gordon (editor)
{\em Heegaard splittings of 3œôòómanifolds (Haifa 2005), Problems}
Geom. Topol. Monographs {\bf 12} (2007) 401-œôòó411.


\bibitem{CG} A.~Casson and C.~Gordon, {\em Manifolds with irreducible Heegaard splittings of arbitrarily
    high genus}, preprint.


\bibitem{FN} A.~Felikson, S.~Natanzon, {\em Double pants decompositions of 2-surfaces},  arXiv:1005.0073v2.

\bibitem{He} J.~Hempel, {\em 3-manifolds as viewed from the curve complex}, Topology {\bf 40}, No. 3 (2001), 631--657.


\bibitem{H1} A.~Hatcher, {\em Pants decompositions of surfaces}, \\  
\url{http://www.math.cornell.edu/~hatcher/Papers/pantsdecomp.pdf} 

\bibitem{HT} A.~Hatcher, W.~Thurston, {\em A presentation for the mapping class group of a closed orientable surface},
Topology {\bf 19} (1980), 221--237. 

\bibitem{J} 
J.~Johnson, {\em Heegaard splittings and the pants complex}. Algebr. Geom. Topol., {\bf 6} (2006), 853--874.

\bibitem{L}
F.~Luo, {\em On Heegaard diagrams}, Math. Res. Lett., {\bf 4} (1997), no 2-3, 365--373.


\end{thebibliography}
\end{document}